\documentclass[12pt]{article}

\usepackage{amsmath, epsfig, cite}
\usepackage{amssymb}
\usepackage{amsfonts}
\usepackage{latexsym}
\usepackage{float}
\usepackage{color}
\usepackage{url}

\newtheorem{thm}{Theorem}[section]

\newtheorem{conj}[thm]{Conjecture}

\numberwithin{equation}{section}

\begin{document}

\begin{center}
{\Large\bf On the divisibility of sums of even powers of\\[10pt]
 $q$-binomial coefficients}
\end{center}

\vskip 2mm \centerline{Ji-Cai Liu and Xue-Ting Jiang}
\begin{center}
{\footnotesize Department of Mathematics, Wenzhou University, Wenzhou 325035, PR China\\
{\tt jcliu2016@gmail.com, xtjiang2021@foxmail.com} }
\end{center}


\vskip 0.7cm \noindent{\bf Abstract.}
We prove the divisibility conjecture on sums of even powers of $q$-binomial coefficients, which was recently proposed by Guo, Schlosser and Zudilin. Our proof relies on two $q$-harmonic series congruences due to Shi and Pan.


\vskip 2mm
\noindent{\it MR Subject Classifications}: 11A07, 11B65, 05A10

\section{Introduction}
In 2012, Guo and Zeng \cite[(1.9)]{gz-ijnt-2012} established the following divisibility result
to prove two amazing supercongruences involving the Ap\'ery numbers:
\begin{align}
\sum_{k=0}^{n-1}{n+k\choose k}^2{n-1\choose k}^2\equiv 0\pmod{n}.\label{aa-2}
\end{align}
Guo and Zeng \cite[Conjecture 5.13]{gz-ijnt-2012} also conjectured a
$q$-analogue of \eqref{aa-2} as follows:
\begin{align}
\sum_{k=0}^{n-1}q^{(n-k)^2}{n+k\brack k}^2{n-1\brack k}^2\equiv q^{(n-1)^2}[n]\pmod{[p]^2_{q^{n/p}}},\label{aa-3}
\end{align}
where $n$ is any power of a prime $p$.

To understand \eqref{aa-3}, we shall recall some notation.
The Gaussian $q$-binomial coefficients are defined as
\begin{align*}
{n\brack k}={n\brack k}_q
=\begin{cases}
\displaystyle\frac{(q;q)_n}{(q;q)_k(q;q)_{n-k}} &\text{if $0\leqslant k\leqslant n$},\\[10pt]
0 &\text{otherwise,}
\end{cases}
\end{align*}
where the $q$-shifted factorials are given by $(a;q)_n=(1-a)(1-aq)\cdots(1-aq^{n-1})$ for $n\ge 1$ and $(a;q)_0=1$. The $q$-integers are defined by $[n]=[n]_q=(1-q^n)/(1-q)$. The $n$th cyclotomic polynomial is given by
\begin{align*}
\Phi_n(q)=\prod_{\substack{1\le k \le n\\[3pt](n,k)=1}}
(q-\zeta^k),
\end{align*}
where $\zeta$ denotes a primitive $n$th root of unity.

For polynomials $A_1(q), A_2(q),P(q)\in \mathbb{Z}[q]$, the $q$-congruence $$A_1(q)/A_2(q)\equiv 0\pmod{P(q)}$$ is understood as $A_1(q)$ is divisible by $P(q)$ and $A_2(q)$ is coprime with $P(q)$. In general, for rational functions $A(q),B(q)\in \mathbb{Q}(q)$ and polynomial $P(q)\in\mathbb{Z}[q]$,
\begin{align*}
A(q)\equiv B(q)\pmod{P(q)}\Longleftrightarrow
A(q)-B(q)\equiv 0\pmod{P(q)}.
\end{align*}

In the past few years, congruences for $q$-binomial coefficients as well as basic hypergeometric series attracted many experts' attention (see, for example, \cite{bachraoui-rj-2021,goro-ijnt-2019,guo-racsam-2020,guo-racsam-2021,gl-jdea-2018,gs-racsam-2021,
gz-am-2019,lw-racsam-2020,liu-crm-2020,lh-bams-2020,lp-aam-2020,wh-pmh-2021,wy-racsam-2021,wei-jcta-2021}).

In 2020, Gu and Guo \cite[Theorem 1.2]{gg-cmj-2020} confirmed the conjectural $q$-congruence \eqref{aa-3} by establishing the following more general result:
\begin{align}
\sum_{k=0}^{n-1}q^{(n-k)^2}{n+k\brack k}^2{n-1\brack k}^2\equiv q[n]\pmod{\Phi_n(q)^2},\label{aa-4}
\end{align}
where $n$ is any positive integer. They also conjectured that \eqref{aa-4} holds modulo $[n]\Phi_n(q)^2$ (see \cite[Conjecture 1.3]{gg-cmj-2020}), which is still open up to now.

Recently, Guo, Schlosser and Zudilin \cite[Theorem 8]{gsz-2021} generalize \eqref{aa-4} in a different form:
\begin{align}
\sum_{k=0}^{n-1}q^{r(n-k)^2+(r-1)k}{n+k\brack k}^{2r}{n-1\brack k}^{2r}
\equiv q[n]\pmod{\Phi_n(q)^2},\label{aa-5}
\end{align}
where $r$ is any positive integer. They pointed out that
via the same method as in the proof of \cite[Theorem 5.3]{gz-ijnt-2012}, one can easily show that
\begin{align}
\sum_{k=0}^{n-1}q^{r(n-k)^2+(r-1)k}{n+k\brack k}^{2r}{n-1\brack k}^{2r}
\equiv 0\pmod{[n]}.\label{a-1}
\end{align}
It follows that \eqref{aa-5} is also true modulo $[n]\Phi_n(q)$. Guo, Schlosser and Zudilin \cite[Conjecture 1]{gsz-2021} also conjectured that both \eqref{aa-5} and the conjecture \cite[Conjecture 1.3]{gg-cmj-2020}) due to Gu and Guo possess the following unified generalization:

\begin{conj}[Guo--Schlosser--Zudilin]\label{gsz}
Let $n$ be a positive integer and $r$ an arbitrary integer. Then
\begin{align}
&\sum_{k=0}^{n-1}q^{r(n-k)^2+(r-1)k}{n+k\brack k}^{2r}{n-1\brack k}^{2r}\notag\\[10pt]
&\equiv q^{(r-1)n+1}[n]-\frac{r(2r-1)(n-1)^2q(1-q)^2}{4}[n]^3\pmod{[n]\Phi_n(q)^3}.\label{conj-1}
\end{align}
\end{conj}

At the end of their paper, Guo, Schlosser and Zudilin \cite{gsz-2021} mentioned that Conjecture \ref{gsz} is even still open modulo $[n]\Phi_n(q)^2$. The motivation of this paper is to give a positive answer to Conjecture \ref{gsz}.

\begin{thm}\label{t-1}
The $q$-congruence \eqref{conj-1} is true.
\end{thm}

The important ingredients in the proof of Theorem \ref{t-1} are the following two $q$-harmonic series congruences:
\begin{align}
\sum_{j=1}^{n-1}\frac{1}{1-q^j}\equiv \frac{n-1}{2}+\frac{(n^2-1)(1-q^n)}{24}\pmod{\Phi_n(q)^2},\label{b-4}
\end{align}
and
\begin{align}
\sum_{j=1}^{n-1}\frac{1}{(1-q^j)^2}\equiv -\frac{(n-1)(n-5)}{12}\pmod{\Phi_n(q)}.\label{b-5}
\end{align}
As Pan \cite[Proposition 2.1]{pan-rmjm-2008} mentioned, the $q$-congruence \eqref{b-4} can be easily proved by the same method as in the proof of \cite[Theorem 1]{sp-amm-2007}. In fact, one can also prove
\eqref{b-5} in a similar way as in the proof of \cite[(5)]{sp-amm-2007}.
In the next section, we shall present the detailed proof of Theorem \ref{t-1}.

\section{Proof of Theorem \ref{t-1}}
Noting that
\begin{align*}
(1-q^{n-j})(1-q^{n+j})=(1-q^n)^2-(1-q^j)^2q^{n-j},
\end{align*}
we obtain
\begin{align}
{n+k\brack k}{n-1\brack k}&=\frac{1}{(q;q)_k^2}\prod_{j=1}^k(1-q^{n-j})(1-q^{n+j})\notag\\[10pt]
&=\frac{1}{(q;q)_k^2}\prod_{j=1}^k\left((1-q^n)^2-(1-q^j)^2q^{n-j}\right)\notag\\[10pt]
&\equiv (-1)^k q^{nk-\frac{k(k+1)}{2}}\left(1-\frac{(1-q^n)^2}{q^n}\sum_{j=1}^k\frac{q^j}{(1-q^j)^2}\right)
\pmod{\Phi_n(q)^4},\label{new-4}
\end{align}
where we have used the fact that $1-q^n\equiv 0\pmod{\Phi_n(q)}$.

By \eqref{new-4} and the binomial theorem, we have
\begin{align*}
{n+k\brack k}^{2r}{n-1\brack k}^{2r}
\equiv q^{2rnk-rk(k+1)}\left(1-\frac{2r(1-q^n)^2}{q^n}\sum_{j=1}^k\frac{q^j}{(1-q^j)^2}\right)\pmod{\Phi_n(q)^4}.
\end{align*}
It follows that
\begin{align}
&\sum_{k=0}^{n-1}q^{r(n-k)^2+(r-1)k}{n+k\brack k}^{2r}{n-1\brack k}^{2r}\notag\\[10pt]
&\equiv q^{n^2r}\sum_{k=0}^{n-1}q^{-k}-2rq^{n(nr-1)}(1-q^n)^2\sum_{k=0}^{n-1}q^{-k}
\sum_{j=1}^k\frac{q^j}{(1-q^j)^2}\pmod{\Phi_n(q)^4}.\label{b-1}
\end{align}

It is clear that
\begin{align}
q^{n^2r}\sum_{k=0}^{n-1}q^{-k}=q^{n(nr-1)+1}[n].\label{b-2}
\end{align}
Furthermore, we have
\begin{align}
\sum_{k=0}^{n-1}q^{-k}
\sum_{j=1}^k\frac{q^j}{(1-q^j)^2}&=\sum_{j=1}^{n-1}\frac{q^j}{(1-q^j)^2}\sum_{k=j}^{n-1}q^{-k}\notag\\[10pt]
&=\frac{1}{q^{n-1}(q-1)}\sum_{j=1}^{n-1}\frac{q^n-q^j}{(1-q^j)^2}\notag\\[10pt]
&=\frac{1}{q^{n-1}(q-1)}\left(\sum_{j=1}^{n-1}\frac{1}{1-q^j}-(1-q^n)\sum_{j=1}^{n-1}\frac{1}{(1-q^j)^2}\right).
\label{b-3}
\end{align}
Substituting \eqref{b-4} and \eqref{b-5} into the right-hand side of \eqref{b-3}, we obtain
\begin{align}
&\sum_{k=0}^{n-1}q^{-k}\sum_{j=1}^k\frac{q^j}{(1-q^j)^2}\notag\\[10pt]
&\equiv
\frac{1}{q^{n-1}(q-1)}\left(\frac{n-1}{2}+\frac{(n-1)(n-3)(1-q^n)}{8}\right)\pmod{\Phi_n(q)^2}.
\label{b-6}
\end{align}
Combining \eqref{b-1}, \eqref{b-2} and \eqref{b-6} gives
\begin{align}
\sum_{k=0}^{n-1}q^{r(n-k)^2+(r-1)k}{n+k\brack k}^{2r}{n-1\brack k}^{2r}
&\equiv q^{n(nr-1)+1}[n]+2rq^{n(nr-2)+1}(1-q)[n]^2\notag\\[10pt]
&\hskip -30mm\times\left(\frac{n-1}{2}
+\frac{(n-1)(n-3)(1-q)[n]}{8}\right)\pmod{\Phi_n(q)^4}.\label{new-3}
\end{align}

Note that
\begin{align*}
q^{sn}&=1-(1-q^{sn})\notag\\[10pt]
&=1-(1-q^n)(1+q^n+q^{2n}+\cdots+q^{(s-1)n})\notag\\[10pt]
&\equiv 1-s(1-q^n)\pmod{\Phi_n(q)^2}.
\end{align*}
It follows that
\begin{align}
q^{sn}&=1-(1-q^n)(1+q^n+q^{2n}+\cdots+q^{(s-1)n})\notag\\[10pt]
&\equiv 1-(1-q^n)\left(s-\frac{s(s-1)(1-q^n)}{2}\right)\pmod{\Phi_n(q)^3}\notag\\[10pt]
&=1-s(1-q)[n]+\frac{s(s-1)(1-q)^2}{2}[n]^2.\label{new-2}
\end{align}
Applying \eqref{new-2} to the right-hand side of \eqref{new-3} yields
\begin{align}
\sum_{k=0}^{n-1}q^{r(n-k)^2+(r-1)k}{n+k\brack k}^{2r}{n-1\brack k}^{2r}
&\equiv q[n]-q(1-q)(r-1)[n]^2\notag\\[10pt]
&\hskip-60mm-\frac{q(1-q)^2(2{n}^{2}{r}^{2}-{n}^{2}r-4n{r}^{2}+2nr+5r-4)}{4}[n]^3\pmod{\Phi_n(q)^4}.
\label{b-9}
\end{align}

On the other hand, by \eqref{new-2} we have
\begin{align}
q^{(r-1)n+1}[n]-\frac{r(2r-1)(n-1)^2q(1-q)^2}{4}[n]^3
&\equiv q[n]-q(1-q)(r-1)[n]^2\notag\\[10pt]
&\hskip-70mm-\frac{q(1-q)^2(2{n}^{2}{r}^{2}-{n}^{2}r-4n{r}^{2}+2nr+5r-4)}{4}[n]^3\pmod{\Phi_n(q)^4}.
\label{b-10}
\end{align}
Finally, comparing the right-hand sides of \eqref{b-9} and \eqref{b-10} gives
\begin{align}
&\sum_{k=0}^{n-1}q^{r(n-k)^2+(r-1)k}{n+k\brack k}^{2r}{n-1\brack k}^{2r}\notag\\[10pt]
&\equiv q^{(r-1)n+1}[n]-\frac{r(2r-1)(n-1)^2q(1-q)^2}{4}[n]^3\pmod{\Phi_n(q)^4}.\label{b-11}
\end{align}
Then the proof of \eqref{conj-1} follows from \eqref{a-1} and \eqref{b-11}.

\vskip 5mm \noindent{\bf Acknowledgments.}
The first author was supported by the National Natural Science Foundation of China (grant 12171370).

\end{document}